\documentclass[notitlepage,a4,10.5pt]{article}

\usepackage{graphicx}

\usepackage{amsmath}%
\usepackage{amsfonts}%
\usepackage{amssymb}%
\usepackage{mathrsfs}
\usepackage{amsthm}
\usepackage{float}

\usepackage{authblk}

\usepackage[left=2.5cm,right=2.5cm,top=3cm,bottom=3cm]{geometry} 

\newcommand{\mc}{\mathcal}
\newcommand{\cp}{\times}

\newcommand{\di}{\nabla\cdot}
\newcommand{\cu}{\nabla\times}

\newcommand{\bol}{\boldsymbol}

\newcommand{\abs}[1]{\left\lvert{#1}\right\rvert}

\newcommand{\lr}[1]{\left({#1}\right)}

\newcommand{\mf}{\mathfrak}
\newcommand{\p}{\partial}

\newcommand{\ti}[1]{\textit{#1}}
\newcommand{\tb}[1]{\textbf{#1}}
\newcommand{\ov}[1]{\mkern 1.5mu\overline{\mkern-1.5mu#1\mkern-1.5mu}\mkern 1.5mu}

\DeclareMathOperator{\arcsinh}{arcsinh}

\begin{document}

\title{Quasisymmetric Magnetic Fields in Asymmetric Toroidal Domains}
\author[1]{Naoki Sato} \author[2]{Zhisong Qu} \author[3]{David Pfefferl\'e} \author[2]{Robert L. Dewar} 
\affil[1]{Graduate School of Frontier Sciences, \protect\\ The University of Tokyo, Kashiwa, Chiba 277-8561, Japan \protect\\ Email: sato\_naoki@edu.k.u-tokyo.ac.jp}
\affil[2]{Mathematical Sciences Institute, The Australian National University, Canberra, ACT 2601, Australia}
\affil[3]{The University of Western Australia, 35 Stirling Highway, Crawley WA 6009, Australia}
\date{\today}
\setcounter{Maxaffil}{0}
\renewcommand\Affilfont{\itshape\small}

    \maketitle
    \begin{abstract}
    We explore the existence of quasisymmetric magnetic fields in asymmetric toroidal domains.  
		These vector fields can be identified with a class of magnetohydrodynamic equilibria in the presence of pressure anisotropy.
		First, using Clebsch potentials, 
		we derive a system of two coupled nonlinear first order partial differential equations 
		expressing a family of quasisymmetric magnetic fields in bounded domains.
		In regions where flux surfaces and surfaces of constant field strength are not tangential, this system can be further reduced to a single degenerate nonlinear second order partial differential equation with externally assigned initial data. 
		Then, we exhibit regular quasisymmetric vector fields which correspond to local solutions of anisotropic magnetohydrodynamics
		in asymmetric toroidal domains such that tangential boundary conditions are fulfilled on a portion of the bounding surface. 
		The problems of boundary shape and locality are also discussed. 
		We find that symmetric magnetic fields can be fitted into asymmetric domains, 
		and that the mathematical difficulty encountered in the derivation of global quasisymmetric magnetic fields    
		lies in the topological obstruction toward global extension affecting local solutions of the governing nonlinear first order partial differential equations.
		\end{abstract}

\section{Introduction}
In the context of ideal magnetohydrodynamics with isotropic pressure, 
steady magnetic confinement of a plasma 
is achieved by balance between Lorentz force and pressure gradient,
\begin{subequations}\label{IMHD}
\begin{align}
&\lr{\cu\bol{B}}\cp\bol{B}=\nabla P,~~~~\nabla\cdot\bol{B}=0~~~~{\rm in}~~\Omega,\label{IMHD1}\\
&\bol{B}\cdot\bol{n}=0~~~~{\rm on}~~\p\Omega.
\end{align}
\end{subequations}
In the notation above, $\bol{B}$ is the magnetic field,  
$P$ the pressure, $\Omega\subset\mathbb{R}^3$ a smooth bounded domain, 
and $\bol{n}$ the unit outward normal to the bounding surface $\p\Omega$. 
As a system of nonlinear first order partial differential equations 
for the Cartesian components $B_x$, $B_y$, $B_z$ and the pressure $P$, 
equation \eqref{IMHD} is twice elliptic and twice hyperbolic \cite{YosYam}. 
Such mixed behavior makes \eqref{IMHD} one of the hardest equations in mathematical physics.
Indeed, while in the presence of a continuous Euclidean symmetry 
the hyperbolic part can be removed to obtain 
an elliptic nonlinear second order partial differential equation for the flux function 
(the Grad-Shafranov equation \cite{Grad58,Eden1,Eden2}),  
the existence of regular solutions of \eqref{IMHD} without such isometries remains an unsolved
problem in the theory of partial differential equations \cite{SatoPPCF}. 

In the following, we shall always use the word symmetry to refer to 
continuous transformations of three dimensional Euclidean space that preserve the Euclidean distance
between points, i.e. superpositions of translations and rotations. 
Asymmetry will then imply absence of such Euclidean isometries. 
Furthermore, a certain equilibrium will be deemed symmetric as long as the magnetic field is, without any requirements on other fields or boundary shape. 
It should be noted that a symmetric solution $\lr{\bol{B},P}$ of \eqref{IMHD} 
does not necessarily imply a symmetric boundary $\p\Omega$.  
For example, a symmetric magnetic field satisfying \eqref{IMHD} can be embedded within an asymmetric bounded domain. 
A symmetric solution always implies a symmetric boundary when the boundary corresponds to an isobaric surface. 
 On these points, see \cite{Sato2020}.
According to a conjecture due to H. Grad, 
well-behaved solutions of \eqref{IMHD} with non-constant pressure in $\Omega$
and constant pressure on $\p\Omega$ should possess a high degree of symmetry \cite{Grad67}. 
Intuitively, this is because the boundary condition for the pressure 
effectively forces magnetic field $\bol{B}$ and current $\cu\bol{B}$ on toroidal surfaces
(hairy ball theorem \cite{Hairy}). 
In the absence of a boundary symmetry, this topological obstruction 
cannot be trivially overcome when one tries to lie $\bol{B}$ and $\cu\bol{B}$ consistently around the torus.
A workaround for this problem consists in 
expanding the class of admissible solutions so that controlled discontinuities are allowed.
In this construction, the domain $\Omega$ is partitioned into a given number of subdomains.
Within each subdomain a constant pressure is assumed, and the corresponding magnetic field
is a Beltrami field (an eigenstate of the curl operator, see \cite{YosGiga,Bruno,Enciso} for the
existence of such solutions), while pressure jumps occur at the boundaries among subdomains.  
This approach has proven effective for the modeling of three-dimensional
magnetohydrodynamic equilibria aimed at stellarator design \cite{Dew,Hud,Qu}  
(stellarators are candidate magnetic confinement devices for nuclear fusion applications that 
do not exhibit boundary symmetry). 

Compared to an axially symmetric tokamak, 
a stellarator offers the advantage that 
the field line twist required to trap charged particles
is not induced by currents within the plasma,
but sustained through asymmetric coils, 
thus enabling improved steady state operation of the machine.
This usually comes at the price of deteriorated confinement, 
and a nontrivial coil design that should be compatible with  
a quasisymmetric confining magnetic field \cite{Helander}, i.e. a
magnetic field $\bol{B}$ whose strength $B=\sqrt{\bol{B}\cdot\bol{B}}$ 
is independent of a direction $\bol{u}$ in space, $\bol{u}\cdot\nabla B=0$ (see below for a rigorous definition of quasisymmetric vector field).  
Here, quasisymmetry is required because it guarantees the existence of a conserved momentum
arising from the invariance along $\bol{u}$ of the guiding center Hamiltonian, 
which depends on the magnetic field only through the strength $B$ \cite{Cary}.
The conserved quantity is expected to enhance particle confinement along the magnetic field.

A quasisymmetric magnetic field is mathematically characterized by the property that
there exists a solenoidal vector field $\bol{u}$ such that
both $\bol{B}$ and the field strength $B$ are Lie-invariant
along $\bol{u}$, i.e.
\begin{equation}
\mf{L}_{\bol{u}}\bol{B}=\bol{0},~~~~\mf{L}_{\bol{u}}B=0,~~~~\mf{L}_{\bol{u}}dV=0,\label{QS1}
\end{equation}
where $\mf{L}$ denotes the Lie derivative and $dV=dxdydz$ the volume element in $\mathbb{R}^3$.
Clearly, a symmetric magnetic field is also quasisymmetric. In this case $\bol{u}=\bol{a}+\bol{b}\cp\bol{x}$ 
is the generator of continuous Euclidean isometries, with $\bol{a},\bol{b}\in\mathbb{R}^3$ constant vectors 
and $\bol{x}$ the position vector in $\mathbb{R}^3$. 
Note that \eqref{IMHD1} can also be expressed through the Lie derivative as $\mf{L}_{\cu{\bol{B}}}\bol{B}=\bol{0}$ and $\mf{L}_{\bol{B}}dV=0$.
The quasisymmetry condition \eqref{QS1} is usually written in vector notation as \cite{Rodriguez2020},
\begin{equation}
\bol{B}\cp\bol{u}=\nabla\zeta,~~~~\bol{u}\cdot\nabla B=0,~~~~\di\bol{u}=0.\label{QS2}
\end{equation}
Here, $\zeta$ is any function. 
Nonetheless, in stellarator applications it is customary to identify $\zeta$ with the flux function $\Psi$ taking constant values on the bounding surface, 
so that both the magnetic field $\bol{B}$ and the quasisymmetry $\bol{u}$ lie on flux surfaces, and the conserved momentum
is a combination of the poloidal and toroidal momenta. 
Unfortunately, the analysis of \cite{Garren} suggests that fulfilling \eqref{IMHD} and \eqref{QS2} simultaneously with a constant pressure at the bounding surface 
in an asymmetric configuration should not be possible (see also \cite{Sengupta,Plunk,Constantin2021JPP,Constantin2021CMP}). 
This is because the Lie-invariance required for a quasisymmetric solution to exist   
cannot be satisfied through the available geometrical degrees of freedom.
A possible strategy to overcome this difficulty is to 
modify the force balance equation (the first equation in \eqref{IMHD1}) by introducing 
pressure anisotropy. In this regard, it has been suggested that pressure anisotropy
could remove the problem of overdetermination encountered in the isotropic setting \cite{Rodriguez2020,Rodriguez2021}.  
In the anisotropic case, anisotropic magnetohydrodynamic equilibria  
are described by the boundary value problem
\begin{equation}
\begin{split}
&\lr{\cu\bol{B}}\cp\bol{B}=\nabla\cdot\Pi,~~~~\di\bol{B}=0~~~~{\rm in}~~\Omega,\\
&\bol{B}\cdot\bol{n}=0~~~~{\rm on}~~\p\Omega,\label{aMHD}
\end{split}
\end{equation}
where in Cartesian coordinates $\lr{x^1,x^2,x^3}=\lr{x,y,z}$ 
the twice contravariant symmetric (here, symmetric does not refer to continuous Euclidean isometries but to the property that the transpose tensor equals the tensor itself) pressure tensor $\Pi$ has components
\begin{equation}
\Pi^{ij}=P_{\perp}\delta^{ij}+\frac{\lr{P_{\parallel}-P_{\perp}}}{B^2}B^iB^j,~~~~i,j=1,2,3.\label{Pi}
\end{equation}
The functions $P_{\perp}$ and $P_{\parallel}$ are  
pressure fields associated with the perpendicular and parallel directions to the magnetic field $\bol{B}$. 
Indeed, setting $b^i=B^i/B$, $i=1,2,3$, one has
\begin{equation}
\di\Pi=\bol{b}\cp\lr{\nabla P_{\perp}\cp\bol{b}}+\lr{\bol{b}\cdot\nabla P_{\parallel}}\bol{b}+\lr{P_{\parallel}-P_{\perp}}\nabla\cdot\lr{\bol{b}\bol{b}}.
\end{equation}  
Hence, the gradient of $P_{\perp}$ generates the pressure force in the direction perpendicular to $\bol{B}$, 
while the gradient of $P_{\parallel}$ generates the pressure force in the direction parallel to $\bol{B}$.
Isotropic magnetohydrodynamic equilibria can be recovered when $P_{\perp}=P_{\parallel}=P$, with $P$ the scalar pressure. 
We also remark that the form of the pressure tensor \eqref{Pi} is rooted in kinetic theory: the pressure fields $P_{\perp}$ and $P_{\parallel}$ represent the averaged squared deviations of the charged particle velocities in the guiding center approximation across and along the magnetic field with respect to the mean perpendicular and parallel flows  \cite{Grad67,Grad66,Dobrott,Iacono}.

Our aim in this paper is to study the boundary value problem \eqref{aMHD} with the quasisymmetry condition
\eqref{QS2} for the magnetic field. It will be shown that, by appropriately choosing the values of 
$P_{\perp}$ and $P_{\parallel}$, it is possible to construct regular quasisymmetric magnetic fields, i.e. asymmetric 
solutions $\bol{B}$ of system \eqref{QS2}, \eqref{aMHD} in a neighborhood $U\subset\Omega$ within an asymmetric toroidal volume $\Omega$ 
such that $\p U\cap\p\Omega\neq\emptyset$.  
We remark that for the time being we shall not be concerned with the feasibility of the obtained equilibria in the laboratory or their stability,
but simply focus on their existence. 

The present paper is organized as follows.
In section 2, general considerations on the geometrical implications of quasisymmetry are given.
In particular, equation \eqref{QS2} is written in terms of a set of Clebsch potentials \cite{YosClebsch}.
This form will be useful to build explicit solutions.
In section 3, we discuss the representation of asymmetric tori in terms of special sets of coordinates. 
In section 4, we derive regular symmetric equilibria that solve system \eqref{aMHD} 
within asymmetric toroidal domains.
In section 5, we construct regular quasisymmetric solutions of system \eqref{aMHD} in a neighborhood $U\subset\Omega$ of an asymmetric toroidal volume $\Omega$ that satisfy tangential boundary conditions on a portion of the boundary $\p U\cap\p\Omega\neq\emptyset$.
Concluding remarks are given in section 5.

\section{General remarks on quasisymmetry}
Consider Cartesian coordinates $\lr{x^1,x^2,x^3}=\lr{x,y,z}$ in $\Omega$. It is convenient to write the pressure fields $P_{\perp}$ and $P_{\parallel}$ as follows:
\begin{equation}
P_{\perp}=\mc{P}-\frac{1}{2}\sigma B^2,~~~~P_{\parallel}=\mc{P}+\frac{1}{2}\sigma B^2.
\end{equation}
Here, the functions $\mc{P}$ and $\sigma$ can be interpreted as a reference pressure field and a deviation (anisotropy) term respectively. 
The Cartesian components of the pressure tensor become
\begin{equation}
\Pi^{ij}=\lr{\mc{P}-\frac{1}{2}\sigma B^2}\delta^{ij}+\sigma B^iB^j,~~~~i,j=1,2,3.
\end{equation}
In the following, we denote with $\p_i$, $i=1,2,3$, the tangent vector in the $x^i$ direction.
Then, using $\nabla\cdot\bol{B}=0$, one can verify the identity
\begin{equation}
\nabla\cdot\Pi=\lr{\p_i\Pi^{ij}}\p_j=\nabla\lr{\mc{P}-\frac{1}{2}\sigma B^2}+\bol{B}\cdot\nabla\lr{\sigma\bol{B}}.
\end{equation}
Hence, force balance now takes the form
\begin{equation}
\lr{1-\sigma}\lr{\nabla\cp\bol{B}}\cp\bol{B}=\nabla\mc{P}-\frac{1}{2}B^2\nabla\sigma+\lr{\bol{B}\cdot\nabla\sigma}\bol{B}.\label{FB}
\end{equation}
In the presence of symmetry, it is known that system \eqref{FB} can be reduced to a well-posed Grad-Shafranov type equation under appropriate ellipticity conditions, among which the requirement  $1-\sigma>0$ 
\cite{Grad67,Iacono,Layden}.
Isotropic equilibria correspond to the case $\mc{P}=P$ and $\sigma=0$. More generally, solutions with constant $\sigma\neq 1$ are qualitatively equivalent to isotropic configurations with $P=\mc{P}/(1-\sigma)$. 
Next, let $P_0\in\mathbb{R}$ be a constant. 
The force balance equation \eqref{FB} can be trivially satisfied by setting $\sigma=1$ and $\mc{P}=P_0/2$ so that
\begin{equation}
P_{\perp}=\frac{1}{2}\lr{P_0-B^2},~~~~P_{\parallel}=\frac{1}{2}\lr{P_0+B^2},\label{PP}
\end{equation}  
With this choice, it follows that finding a solution of \eqref{aMHD} under the quasisymmetry condition \eqref{QS2} now amounts to 
solving the boundary value problem
\begin{equation}
\begin{split}
&\nabla\cdot\bol{B}=0,~~~~\bol{B}\cp\bol{u}=\nabla\zeta,~~~~\bol{u}\cdot\nabla B=0,~~~~\nabla\cdot\bol{u}=0~~~~{\rm in}~~\Omega,\\
&\bol{B}\cdot\bol{n}=0~~~~{\rm on}~~\p\Omega.\label{aMHD2}
\end{split}
\end{equation}
Notice that these are just the equations satisfied by a quasisymmetric vector field tangential to  $\p\Omega$. 
Furthermore, observe that in this construction the pressure fields $P_{\perp}$ and $P_{\parallel}$ 
are not constrained at the boundary. Instead, they are functions of the magnetic energy density.
In particular, the perpendicular pressure field $P_{\perp}$ exactly balances the magnetic pressure since $P_{\perp}+B^2/2=P_0/2$, 
the mechanical pressure along the magnetic field $P_{\parallel}$ satisfies $P_{\parallel}+B^2/2=P_0/2+B^2$, 
while the total energy density associated with pressure and magnetic fields is $P_{\perp}+P_{\parallel}+B^2/2=P_0+B^2/2$.    
If one further demands that, for example, $P_{\perp}=0$ on $\p\Omega$, this results in an additional 
boundary condition for the magnetic field strength, $\bol{B}^2=P_0$ on $\p\Omega$. 
We will briefly return to this point later on, 
although the present study will be mainly focused on the case in which the pressure fields are not constrained at the bounding surface.
We also remark that, if the quasisymmetry condition \eqref{QS2} is not imposed, 
any solenoidal magnetic field satisfying tangential boundary conditions 
represents a solution of anisotropic magnetohydrodynamics \eqref{aMHD} as long as the pressure fields
are chosen as in \eqref{PP}.

It is useful to briefly discuss the physical consistency of equation \eqref{PP}. 
First, observe that both $P_{\perp}$ and $P_{\parallel}$ can be made non-negative
by taking a sufficiently large constant $P_0$. Now suppose that the curvature $\abs{\bol{b}\cdot\nabla\bol{b}}$ 
is a small quantity, and rewrite the force balance equation $\lr{\nabla\cp\bol{B}}\cp\bol{B}=\nabla\cdot\Pi$ in the following form:
\begin{equation}
-\frac{1}{2}\bol{b}\cp\lr{\nabla B^2\cp\bol{b}}+B^2\bol{b}\cdot\nabla\bol{b}=\bol{b}\cp\lr{\nabla P_{\perp}\cp\bol{b}}+\lr{P_{\parallel}-P_{\perp}}\bol{b}\cdot\nabla\bol{b}+\left[\bol{b}\cdot\nabla P_{\parallel}
+\lr{P_{\parallel}-P_{\perp}}\lr{\nabla\cdot\bol{b}}\right]\bol{b}.
\end{equation}  
Collecting small terms involving the curvature $\abs{\bol{b}\cdot\nabla\bol{b}}$, one thus arrives at the system
\begin{subequations}
\begin{align}
&B^2=P_{\parallel}-P_{\perp},\\
&\bol{b}\cp\left[\nabla\lr{P_{\perp}+\frac{B^2}{2}}\cp\bol{b}\right]=\bol{0},\\
&\bol{b}\cdot\nabla P_{\parallel}
+\lr{P_{\parallel}-P_{\perp}}\lr{\nabla\cdot\bol{b}}=0.
\end{align}
\end{subequations}
By eliminating $P_{\parallel}$, these equations reduce to
\begin{subequations}
\begin{align}
&\bol{b}\cp\left[\nabla\lr{P_{\perp}+\frac{B^2}{2}}\cp\bol{b}\right]=\bol{0},\\
&\bol{b}\cdot\nabla\lr{P_{\perp}+\frac{B^2}{2}}=0,
\end{align}
\end{subequations}
which imply equation \eqref{PP} as desired.
Hence, equation \eqref{PP} corresponds to a magnetic configuration in which 
the magnetic field curvature $\abs{\bol{b}\cdot\nabla\bol{b}}$ 
represents a higher-order term in the force balance equation.  
 
Now consider system \eqref{aMHD2}. We must distinguish two cases. 
\begin{enumerate}
\item The vector fields $\bol{B}$ and $\bol{u}$ are linearly dependent. Then, $\nabla\zeta=\bol{0}$,
and the magnetic field $\bol{B}$ is self-quasisymmetric. Indeed, system \eqref{aMHD2} reduces to the 
magnetic differential equation
\begin{equation}\label{SQS}
\begin{split}
&\nabla\cdot\bol{B}=0,~~~~\bol{B}\cdot\nabla B=0~~~~{\rm in}~~\Omega,\\
&\bol{B}\cdot\bol{n}=0~~~~{\rm on}~~\p\Omega. 
\end{split}
\end{equation} 
Note that self-quasisymmetry implies that field strength does not change along field lines. 
Hence, guiding center motion is not accelerated in the direction of the magnetic field,
and the conserved momentum is the momentum parallel to $\bol{B}$.    
Due to the regularity of the boundary $\p\Omega$, the unit outward normal $\bol{n}$ can be expressed as
\begin{equation}
\bol{n}=\frac{\nabla \Psi}{\abs{\nabla \Psi}},
\end{equation}
where $\Psi$ is a single-valued function corresponding to the usual flux function. 
For the magnetic field $\bol{B}$ to possess nested flux surfaces, we demand that
\begin{equation}
\bol{B}\cdot\nabla\Psi=0~~~~{\rm in}~~\Omega.
\end{equation}
or, denoting with $\Theta$ a possibly multivalued (angle) variable, 
\begin{equation}
\bol{B}=\nabla\Psi\cp\nabla\Theta~~~~{\rm in}~~\Omega.\label{BSQS}
\end{equation}
Recall that, due to the solenoidal nature of $\bol{B}$, the 
variable $\Theta$ can always be chosen as a single-valued function in a sufficiently small 
region $U\subset\Omega$ (Lie-Darboux theorem \cite{DeLeon}). 
In the following, we shall refer to functions associated with the representation of vector fields 
such as $\Psi$ and $\Theta$ as Clebsch potentials. 
It should be noted that Clebsch representations of the form \eqref{BSQS} always correspond to helicity-free configurations since the vector potential $\bol{A}$ associated with the magnetic field $\bol{B}=\cu\bol{A}$ can be chosen as an integrable vector field, $\bol{A}=\Psi\nabla\Theta$. 
Using \eqref{BSQS}, system \eqref{SQS} reduces to a single nonlinear first-order 
partial differential equation for the variable $\Theta$, 
\begin{equation}
\abs{\nabla\Psi\cp\nabla\Theta}^2=E\lr{\Psi,\Theta}~~~~{\rm in}~~\Omega,\label{BSQS2}
\end{equation}
where $E=E\lr{\Psi,\Theta}$ is any non-negative single-valued function of the variables $\Psi$ and $\Theta$. 
Assume a non-vanishing magnetic field, $E\neq 0$. Then, by defining
the function
\begin{equation}
\Theta'=\int\frac{d\Theta}{\sqrt{E\lr{\Psi,\Theta}}},
\end{equation}
equation \eqref{BSQS2} can be written as
\begin{equation}
\abs{\nabla\Psi\cp\nabla\Theta'}^2=1~~~~{\rm in}~~\Omega.\label{BSQS3}
\end{equation}
Notice that equation \eqref{BSQS3} resembles the eikonal equation. Indeed, it is equivalent to
\begin{equation}
\abs{\nabla\Theta'}^2=\frac{1}{\abs{\nabla\Psi}^2}+\frac{\lr{\nabla\Psi\cdot\nabla\Theta'}^2}{\abs{\nabla\Psi}^2}~~~~{\rm in}~~\Omega.\label{BSQS4}
\end{equation}
In general, the solution of a first-order partial differential equation such as \eqref{BSQS4} involves
the integration of the associated characteristic system of ordinary differential equations.  
Even if such solution is obtained, it usually has a local nature. 
Therefore, the existence of a global self-quasisymmetric magnetic field in $\Omega$ is 
contingent upon the possibility of extending and/or patching local solutions consistently within $\Omega$. 
We also remark that, if one does not enforce boundary conditions, 
constructing self-quasisymmetric fields becomes a simpler task. 
For example, denoting with $\lr{r,\phi,z}$ cylindrical coordinates, the vector field
\begin{equation}
\bol{B}=f\lr{r,\frac{z}{r}-\phi}\nabla r\cp\nabla\lr{\frac{z}{r}-\phi},\label{BHel}
\end{equation}
with $f$ an arbitrary function of $r$ and $z/r-\phi$, satisfies $\nabla\cdot\bol{B}=0$, $\bol{B}\cdot\nabla B=0$,
and it does not possess continuous Euclidean symmetries in general because the equation $\mf{L}_{\bol{u}}\bol{B}=\bol{0}$ 
with $\bol{u}=\bol{a}+\bol{b}\cp\bol{x}$ does not have solutions such that $\bol{u}\neq\bol{0}$. In particular, the quantity $z/r-\phi$ 
does not correspond to an Euclidean helical symmetry due to the radial dependence.   
In addition, the vector field \eqref{BHel} is tangential to cylindrical surfaces (see figure \ref{fig1}). 
Finally, note that in order to validate the asymmetry of a solution in most cases it is easier to 
verify the violation of the condition $\mf{L}_{\bol{u}}B^2=\bol{u}\cdot\nabla B^2=0$. 
Indeed, choosing a curvilinear coordinate system $\lr{y^1,y^2,y^3}$ such that $\bol{u}=\p/\p y^3$ 
and denoting with $B^k$ and $g_{k\ell}$, $k,\ell=1,2,3$, the components of magnetic field and the 
Euclidean metric tensor with respect to the new coordinates, 
it follows that $\bol{u}\cdot\nabla B^2=2g_{k\ell}B^k\p B^{\ell}/\p y^3$ 
since by hypothesis $\bol{u}$ is an Euclidean isometry such that $\p g_{k\ell}/\p y^3=0$, $k,\ell=1,2,3$.
On the other hand, $\mf{L}_{\bol{u}}\bol{B}=\lr{\p B^i/\p y^3}\lr{\p/\p y^i}$, which implies that 
$\mf{L}_{\bol{u}}\bol{B}$ cannot vanish as long as $\mf{L}_{\bol{u}}B^2$ is different from zero.

\begin{figure}[h]
\hspace*{-0cm}\centering
\includegraphics[scale=0.44]{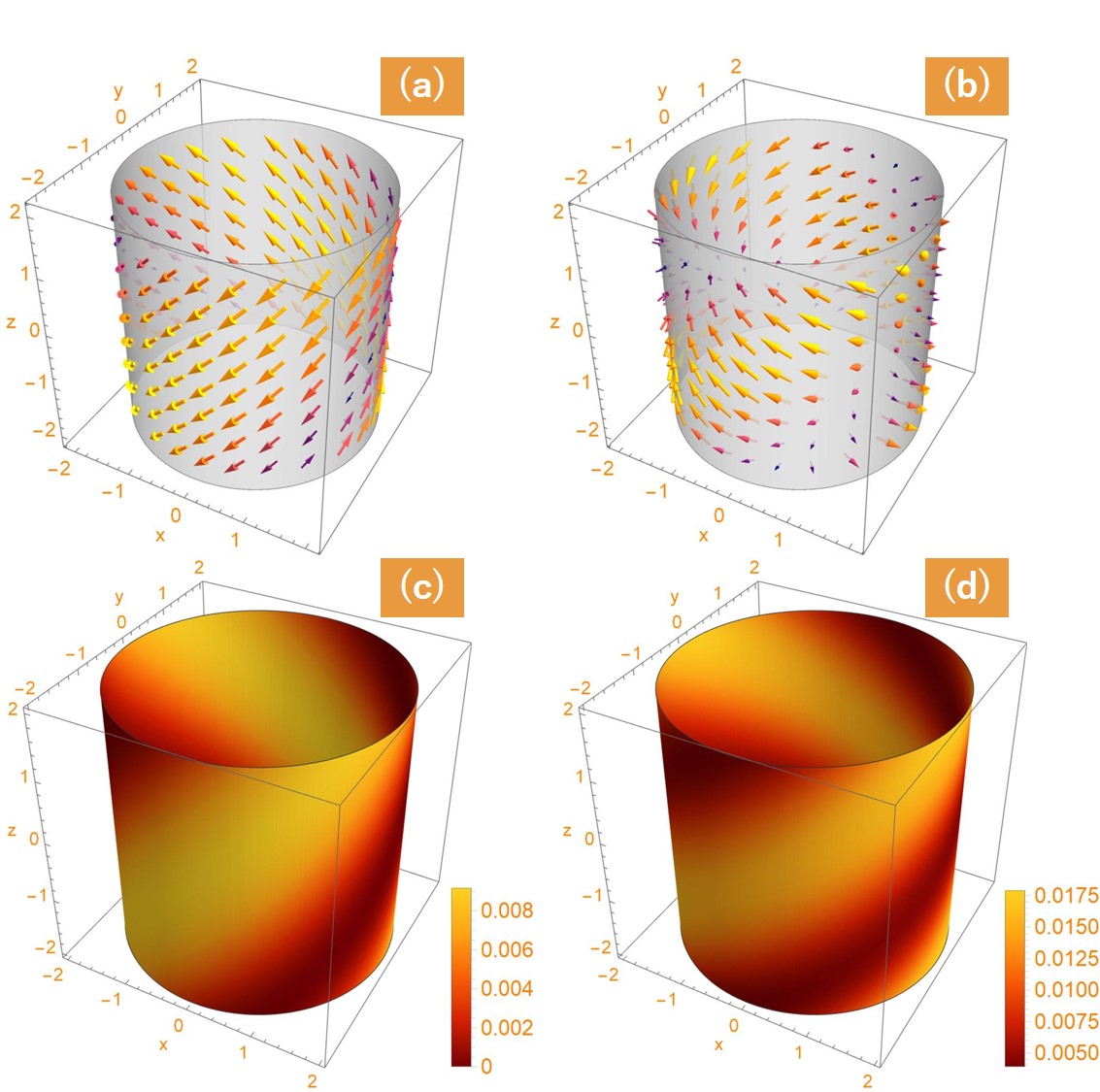}
\caption{\footnotesize (a) Plot of the self-quasisymmetric magnetic field $\bol{B}=e^{-r}\sin\lr{z/r-\phi}\nabla r\cp\nabla\lr{z/r-\phi}$ on a level set of $r$. (b) Plot of the current field $\nabla\cp\bol{B}$. (c) Plot of the magnetic field strength $B^2$. Notice that $B^2$ does not change along $\bol{B}$ (compare with (a)). (d)  Plot of the current field strength $\abs{\cu\bol{B}}^2$.}
\label{fig1}
\end{figure}

\item The vector fields $\bol{B}$ and $\bol{u}$ are linearly independent, implying that $\nabla\zeta\neq\bol{0}$. From the second equation in \eqref{aMHD2}, one sees that both $\bol{B}$ and $\bol{u}$ must be orthogonal to $\nabla\zeta$. 
However, notice that on the boundary $\nabla\zeta$ is not necessarily aligned with the unit outward normal $\bol{n}$. 
Hence, in general, $\zeta$ should not be identified with the flux function $\Psi$ taking constant values on $\p\Omega$.  
Denoting with $\vartheta$ and $\rho$ two possibly multivalued (angle) variables, 
and recalling that $\bol{B}$ and $\bol{u}$ are solenoidal vector fields, 
the space of solutions of system \eqref{aMHD2} can be restricted to 
\begin{equation}
\bol{B}=\nabla\zeta\cp\nabla\vartheta,~~~~\bol{u}=\nabla\zeta\cp\nabla\rho.\label{BU}
\end{equation}
Since $\bol{B}\cdot\nabla\Psi=0$ in $\Omega$, from \eqref{BU} one has 
\begin{equation}
\Psi=\Psi\lr{\zeta,\vartheta}.\label{Psi2}
\end{equation}
For a given flux function $\Psi=\Psi\lr{\bol{x}}$, we assume that equation \eqref{Psi2} can be inverted to obtain
\begin{equation}
\zeta=\zeta\lr{\vartheta,\bol{x}}.\label{psi}
\end{equation} 
Then, upon substituting equations \eqref{BU} and \eqref{psi} into equation \eqref{aMHD2}, 
system \eqref{aMHD2} reduces to the following coupled nonlinear first order partial differential equations for the Clebsch potentials $\vartheta$ and $\rho$
\begin{equation}\label{aMHD3}
\nabla\zeta\cp\nabla\vartheta\cdot\nabla\rho=1,~~~~
\abs{\nabla\zeta\cp\nabla\vartheta}^2=F\lr{\zeta,\rho}~~~~{\rm in}~~\Omega,
\end{equation}
where $F$ is any non-negative single-valued function of the Clebsch potentials $\zeta\lr{\vartheta,\bol{x}}$ and $\rho$.
Notice that, since equation \eqref{aMHD3} is a first-order system, 
the main hurdle toward finding an integral is again represented by the possibility of extending a local solution 
to the whole $\Omega$. 
If Clebsch potentials $\vartheta$ and $\rho$ can be determined so that \eqref{aMHD3} is satisfied, the corresponding
vector fields \eqref{BU} define a quasisymmetric solution of anisotropic magnetohydrodynamics.
Finally, observe that system \eqref{aMHD3} can be further reduced to a single nonlinear second-order partial differential equation in regions $V\subset\Omega$ where $\nabla\zeta\cp\nabla B\neq\bol{0}$. 
Indeed, the condition $\bol{u}\cdot\nabla B=0$ then implies that $\rho=\rho\lr{\zeta,B^2}$ with $B^2=\abs{\nabla\zeta\cp\nabla\theta}^2$.
Therefore, the remaining equation $\bol{B}\cp\bol{u}=\nabla\zeta$ can be written as
\begin{equation}
\nabla\zeta\cp\nabla\theta\cdot\nabla\abs{\nabla\zeta\cp\nabla\theta}^2=\frac{1}{\frac{\p\rho}{\p B^2}}~~~~{\rm in}~~V.\label{GSQS}
\end{equation} 
For a given $\rho=\rho\lr{\zeta,B^2}$, a quasisymmetric solution in $V$ can thus be obtained by solving 
the equation above for the unknown variable $\theta$. 
We remark that equation \eqref{GSQS} is degenerate because it is invariant 
under the transformation $\theta\rightarrow\theta+\xi$, with $\xi=\xi\lr{\zeta}$ an arbitrary function of $\zeta$.
Unfortunately, equation \eqref{GSQS} cannot hold in the whole $\Omega$ (that is we never have $V=\Omega$). 
To see this, consider again system \eqref{aMHD2} and suppose that $\nabla\zeta\cp\nabla B\neq\bol{0}$ in $\Omega$. Since both $\bol{B}$ and $\bol{u}$ lie on surfaces of constant $\zeta$, the conditions $\nabla\cdot\bol{u}=0$ and $\bol{u}\cdot\nabla B=0$ imply that the quasisymmetry must satisfy
\begin{equation}
    \bol{u}=\sigma\lr{\zeta,B}\nabla\zeta\cp\nabla B,
\end{equation}
where $\sigma=\sigma\lr{\zeta,B}$ is some function of $\zeta$ and $B$. 
Using the property $\bol{B}\cdot\nabla\zeta=0$, the equation $\bol{B}\cp\bol{u}=\nabla\zeta$ therefore reduces to
\begin{equation}
\sigma\bol{B}\cdot\nabla B=1.\label{div0}
\end{equation}
Due to the regularity of the magnetic field and its derivatives, we must have $\sigma\neq 0$ for \eqref{div0} to hold. 
Using $\nabla\cdot\bol{B}=0$ and $\bol{B}\cdot\bol{n}=0$ on $\p\Omega$, the equation above implies that
\begin{equation}
    0=\int_{\Omega}\nabla\cdot\lr{B\bol{B}}dV=\int_{\Omega}\frac{dV}{\sigma}.\label{div}
\end{equation}
    For the quasisymmetry $\bol{u}$ to be a continuous function, the fraction $1/\sigma$ must be continuous.
Hence, equation \eqref{div} can be satisfied only if there exists at least one point $\bol{x}^\ast\in\Omega$ such that $1/\sigma\lr{\bol{x}^\ast}=0$, implying $\abs{\bol{u}\lr{\bol{x}^\ast}}=\infty$. This contradicts the regularity of $\bol{u}$. 
Therefore, a global regular quasisymmetric configuration cannot exist in the case $\nabla\zeta\cp\nabla B\neq\bol{0}$. In other words, if a global regular quasisymmetric configuration exists, there must be regions/points where $\nabla\zeta\cp\nabla B=\bol{0}$. Evidently, the property $\rho=\rho\lr{\zeta,B^2}$ breaks down at those points where $\nabla\zeta\cp\nabla B=\bol{0}$, and equation \eqref{GSQS} does not hold there.  


\end{enumerate}

\section{Harmonic orthogonal coordinates and asymmetric tori}
\begin{figure}[h]
\hspace*{-0cm}\centering
\includegraphics[scale=0.46]{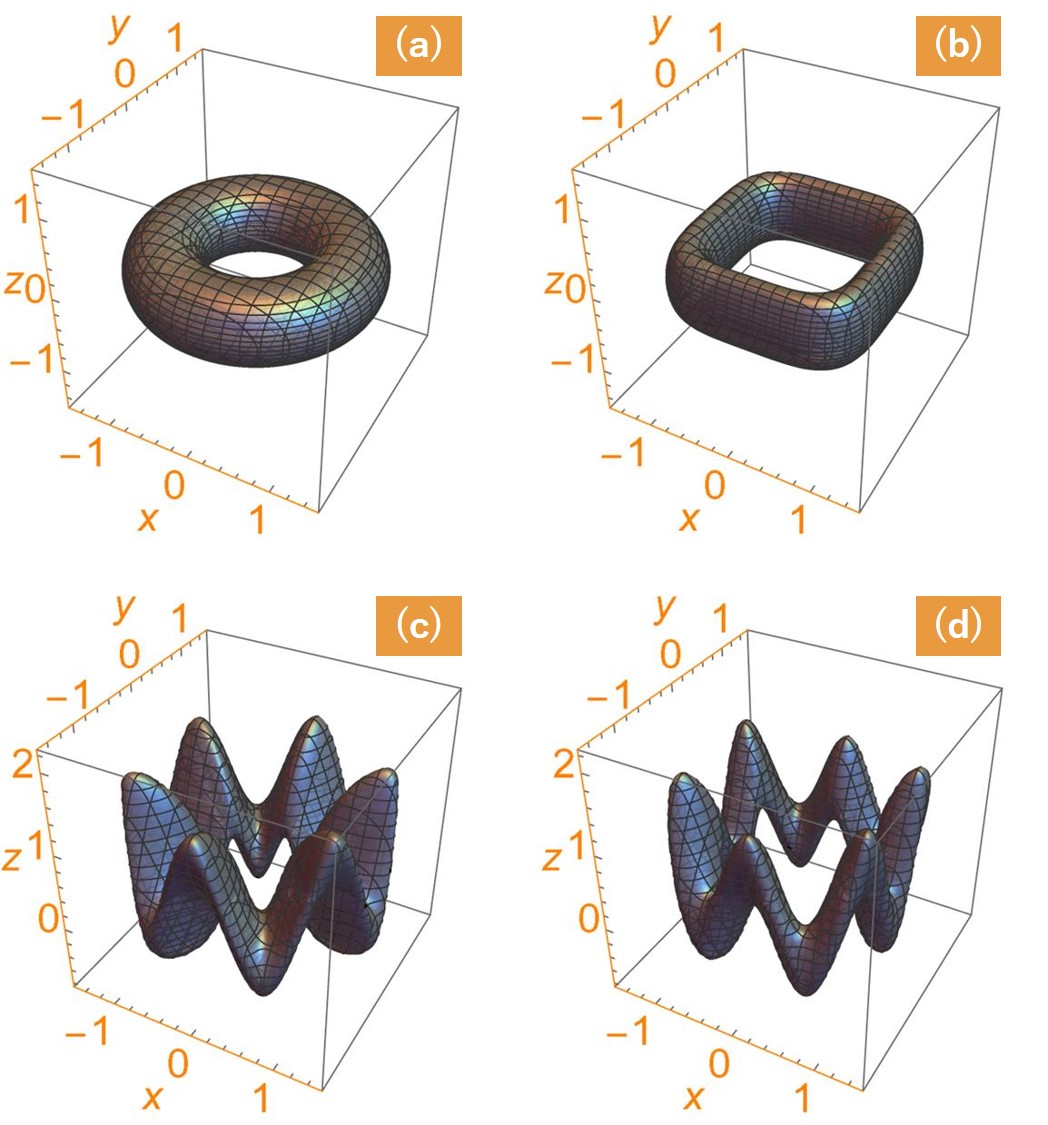}
\caption{\footnotesize (a) Axially symmetric torus $\Psi_{\rm ax}=0.1$ with $r_0=1$. (b) Asymmetric torus $\Psi_{\rm T}=0.1$ with $\mu=\sqrt{x^4+y^4}$, $\mu_0=1$, and $h=0$. (c) Asymmetric torus $\Psi_{\rm T}=0.1$ with $\mu=r$, $\mu_0=1$, and $h=r\sin^2\lr{3\phi}$. (d) Asymmetric torus $\Psi_{\rm T}=0.1$ with $\mu=\sqrt{x^4+y^4}$, $\mu_0=1$, and $h=r\sin^2\lr{3\phi}$.}
\label{fig2}
\end{figure}

In this section, we are concerned with the representation of asymmetric tori 
through special sets of coordinates. Such representation will be used
in the following sections to derive symmetric and quasisymmetric anisotropic magnetohydrodynamic equilibria in asymmetric tori. 
Let $\lr{r,\phi,z}$ denote cylindrical coordinates.
First, consider an axially symmetric torus.  
Such surface is defined as a level set of the function
\begin{equation}
\Psi_{\rm ax}=\frac{1}{2}\left[\lr{r-r_0}^2+z^2\right],\label{axT}
\end{equation} 
with $r_0$ a positive real constant representing the major radius of the torus (the distance from the origin to the center of the torus). 
The axial symmetry is described by the vector field $\bol{u}=\p_{\phi}=r^2\nabla\phi$. Indeed, $\mf{L}_{\p_{\phi}}\Psi_{\rm ax}=0$. 
The axially symmetric torus \eqref{axT} can be thought as a 
special case of a more general family of toroidal surfaces in $\mathbb{R}^3$, which 
correspond to level sets of the function
\begin{equation}
\Psi_{\rm T}=\frac{1}{2}\left[\lr{\mu-\mu_0}^2+\lr{z-h}^2\right].\label{T}
\end{equation}
Here, $\mu=\mu\lr{x,y}$ is a single-valued function in the $\lr{x,y}$ plane measuring the distance from the origin in $\mathbb{R}^2$
and $\mu_0$ a positive real constant corresponding to the major radius of the torus (more generally  
$\mu_0$ can be a single-valued function, although this case is not considered here).  
The more the level sets of $\mu$ differ from circles, the more the torus $\Psi_{\rm T}$ 
departs from axial symmetry in the $\lr{x,y}$ plane.
Similarly, the single-valued function $h=h\lr{\bol{x}}$ expresses the deviation of the toroidal axis from the $\lr{x,y}$ plane. 
Figure \ref{fig2} shows an axially symmetric torus \eqref{axT} and three asymmetric tori \eqref{T} 
obtained from different choices of $\mu$ and $h$. The asymmetry of the tori (b), (c), and (d) in figure \ref{fig1} can be verified by observing that
these surfaces are not Lie-invariant with respect to the generator of continuous Euclidean isometries, that is 
the equation $\mf{L}_{\bol{u}}\Psi_{\rm T}=0$ with $\bol{u}=\bol{a}+\bol{b}\cp\bol{x}$ does not have solution for any nontrivial choice
of the constant vectors $\bol{a},\bol{b}\in\mathbb{R}^3$.
		
If quasisymmetry is not	imposed, solutions of \eqref{aMHD} in asymmetric tori \eqref{T} can be
easily obtained by enforcing \eqref{PP} and setting
\begin{equation}	
\bol{B}=\nabla\Psi_{\rm T}\cp\nabla \Phi,
\end{equation}		
with $\Phi$ an arbitrary function whose gradient field is linearly independent of $\nabla\Psi_{\rm T}$. 
In this case, a sufficient condition for the perpendicular pressure $P_{\perp}=\lr{P_0-B^2}/2$ to be constant on $\p\Omega$ (which corresponds to a level set of $\Psi_{\rm T}$ by construction) is that $\Phi$ is chosen by requiring that $\nabla\Psi_{\rm T}\cp\nabla P_{\perp}=\bol{0}$, or, substituting the expression for $P_{\perp}$,
\begin{equation}
\abs{\nabla\Psi_{\rm T}\cp\nabla\Phi}^2=E\lr{\Psi_{\rm_T}}~~~~{\rm in}~~\Omega,\label{Pb}
\end{equation}
where $E=E\lr{\Psi_{\rm T}}$ is a non-negative function of $\Psi_{\rm T}$. 
Of course, if $P_{\perp}$ is constant on $\p\Omega$, so is $P_{\parallel}$ due to \eqref{PP}.		
Regarding the solvability of \eqref{Pb}, considerations analogous to those pertaining to equation \eqref{BSQS2} apply.
Furthermore, by comparison with equation \eqref{BSQS2}, 
it is clear that a magnetic field fulfilling \eqref{Pb} is also self-quasisymmetric 
since $\bol{B}\cdot\nabla B=0$ (recall equation \eqref{SQS}). 
Solutions of \eqref{Pb} can be obtained easily in axially symmetric tori \eqref{axT}. 
Indeed, first observe that the vector fields
\begin{equation}
\bol{B}=\lambda\nabla\phi,\label{Bax}
\end{equation}
with 
$\lambda=\lambda\lr{r,z}$ are self-quasisymmetric because they satisfy \eqref{SQS} when $\Omega$ is an axially symmetric torus \eqref{axT}. Then, equation \eqref{Pb} holds provided that $\lambda^2=r^2E\lr{\Psi_{\rm ax}}$.
Of course, the self-quasisymmetry of \eqref{Bax} corresponds to the usual rotational symmetry because $\mf{L}_{\p_{\phi}}\bol{B}=\bol{0}$.
However, the construction leading to \eqref{Bax} can be generalized to a certain class of asymmetric tori and, in particular,
to derive symmetric solutions in asymmetric toroidal domains.
To see this, notice that the orthogonal coordinates $\lr{\log{r},\phi,z}$ represent a special case of a larger family of
orthogonal harmonic coordinates $\lr{x^1,x^2,x^3}=\lr{\mu,\nu,z}$ satisfying the following properties in $\Omega$:
\begin{equation}
\nabla x^i\cdot\nabla x^j=\delta^{ij},~~~~\Delta x^i=0,~~~~\abs{\nabla x^1}=\abs{\nabla x^2},~~~~\abs{\nabla x^3}=1,~~~~i,j=1,2,3.\label{HOC}
\end{equation}
The coordinates $\mu$ and $\nu$ can be constructed as follows.  
Let $\Gamma_1$ and $\Gamma_2$ denote two smooth non-intersecting Jordan (simple and closed) curves in $\mathbb{R}^2$.  
Let $\Sigma_1$ and $\Sigma_2$ be the areas enclosed by $\Gamma_1$ and $\Gamma_2$ respectively.
Assume that $\Sigma_2\subset\Sigma_1$ and define $\Gamma=\Gamma_1\cup\Gamma_2$ and $\Sigma=\Sigma_1-\ov{\Sigma}_2$, where the bar denotes the closure of a set. 
Then, the coordinate $\mu=\mu\lr{x,y}$ is obtained as solution of the
following Dirichlet boundary value problem for the Laplace equation, 
\begin{equation}
\Delta\mu=0~~~~{\rm in}~~\Sigma,~~~~\mu=c_1~~~~{\rm on}~~\Gamma_1,~~~~\mu=c_2~~~~{\rm on}~~\Gamma_2,~~~~c_1,c_2\in\mathbb{R},~~~~c_1\neq c_2.
\end{equation}
Observe that $\mu$ can be thought of as a measure of the distance from the origin in $\mathbb{R}^2$.
Next, the coordinate $\nu=\nu\lr{x,y}$ is chosen as the harmonic conjugate of $\mu$.    
In particular, denoting with $\tilde{\bol{n}}$ the unit outward normal to the boundary $\Gamma$,  
the coordinate $\nu$ is a multivalued (angle) variable whose gradient $\nabla\nu$ is a harmonic vector field in $\Sigma$, i.e.
\begin{equation}
\Delta\nu=0~~~~{\rm in}~~\Sigma,~~~~\nabla\nu\cdot\tilde{\bol{n}}=0~~~~{\rm on}~~\Gamma. 
\end{equation}
Finally, equation \eqref{HOC} is satisfied in virtue of the Cauchy-Riemann equations.

\begin{figure}[h]
\hspace*{-0cm}\centering
\includegraphics[scale=0.44]{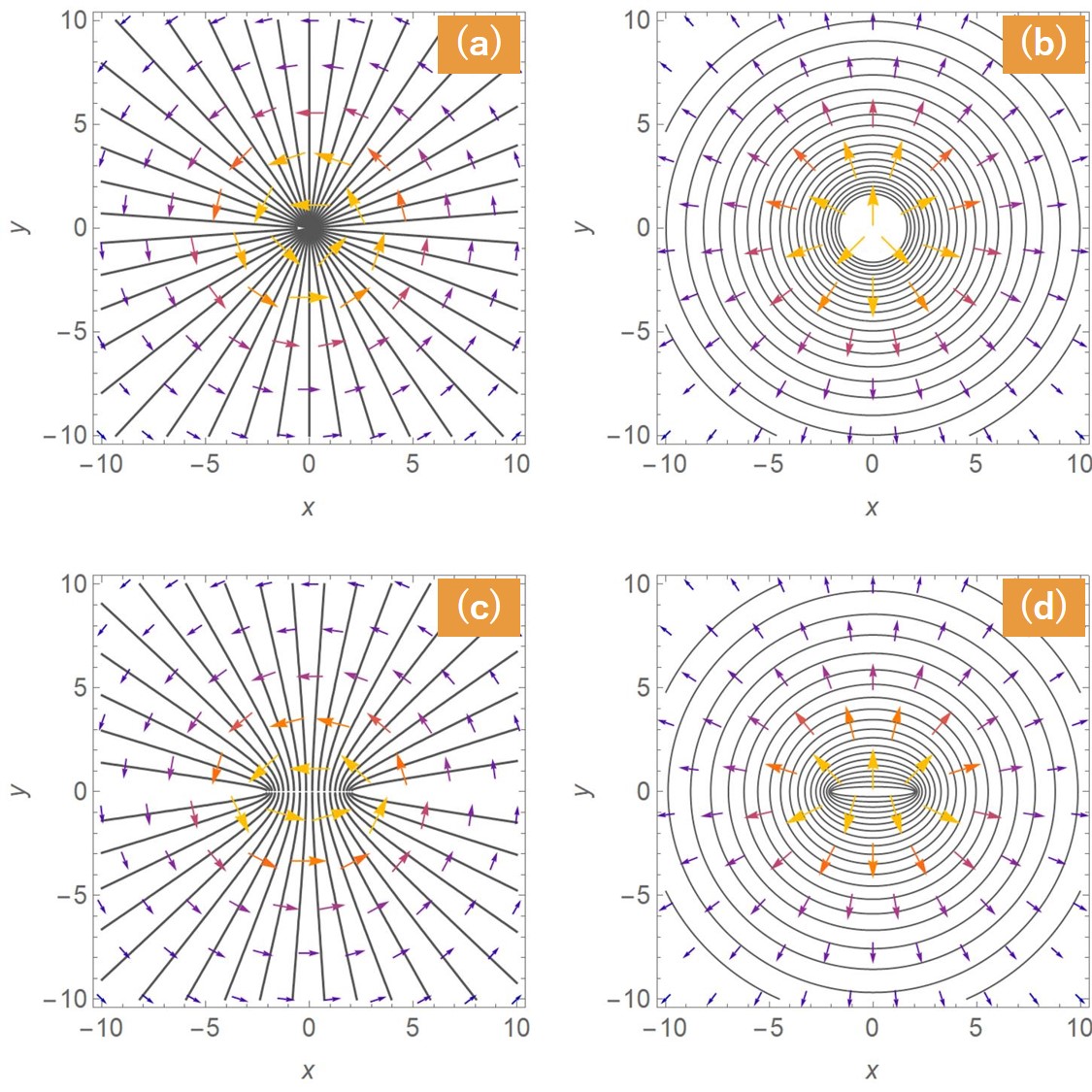}
\caption{\footnotesize (a) Contours of $\phi$ and gradient field $\nabla\phi$. (b) Contours of $\log{r}$ and gradient field $\nabla\log{r}$. (c) Contours of $\nu$ as defined in \eqref{nuEll} and gradient field $\nabla\nu$. (d) Contours of $\mu$ as defined in \eqref{muEll} and gradient field $\nabla\mu$. In these plots, $a=2$.}
\label{fig3}
\end{figure}

As an example, suppose that	$\Gamma_1$ and $\Gamma_2$ are the boundaries of the horizontal sections of two confocal elliptic cylinders 
centered at the origin in $\mathbb{R}^3$. 	
Then, the orthogonal harmonic coordinates $\lr{\mu,\nu,z}$ correspond to the so-called elliptic cylindrical coordinates, which
are related to the Cartesian coordinates $\lr{x,y,z}$ by the transformation
\begin{equation}
x=a\cosh\mu\cos\nu,~~~~y=a\sinh\mu\sin\nu,
\end{equation}
with $a$ a positive real constant such that the foci of the elliptic sections of the cylinders are located at $\bol{x}=\lr{a,0,z}$ and $\bol{x}=\lr{-a,0,z}$.  
Introducing the quantity
\begin{equation}
\delta=\frac{1}{a^2\abs{\nabla\mu}^2}=\frac{1}{a^2\abs{\nabla\nu}^2}=\sin^2\nu+\sinh^2\mu=\sqrt{\lr{1-\frac{x^2+y^2}{a^2}}^2+\frac{4y^2}{a^2}},
\end{equation}
one can derive the following expressions for the quadrant $x,y\geq 0$,
\begin{subequations}\label{EllC}
\begin{align}
&\mu=\arcsinh\left[\frac{1}{\sqrt{2}}\sqrt{-1+\frac{x^2+y^2}{a^2}+\delta}\right],\label{muEll}\\
&\nu=\arcsin\left[\frac{1}{\sqrt{2}}\sqrt{1-\frac{x^2+y^2}{a^2}+\delta}\right],\label{nuEll}
\end{align}
\end{subequations}
and also
\begin{subequations}
\begin{align}
&\nabla\mu=\frac{\sinh\mu\cos\nu\nabla x+\cosh\mu\sin\nu\nabla y}{a\lr{\sin^2\nu+\sinh^2\mu}}=\frac{\frac{\sinh\mu}{\cosh\mu}x\nabla x+\frac{\cosh\mu}{\sinh\mu}y\nabla y}{a^2\lr{\sin^2\nu+\sinh^2\mu}},\\
&\nabla\nu=\frac{\sinh\mu\cos\nu\nabla y-\cosh\mu\sin\nu\nabla x}{a\lr{\sin^2\nu+\sinh^2\mu}}=
\frac{\frac{\sinh\mu}{\cosh\mu} x\nabla y-\frac{\cosh\mu}{\sinh\mu} y\nabla x}{a^2\lr{\sin^2\nu+\sinh^2\mu}}.
\end{align}
\end{subequations}
In the other quadrants, the coordinate $\nu$ can be defined so that  
$\nu\in[0,2\pi)$ when moving around a $\mu$ contour such as $\Gamma_1$ or $\Gamma_2$.
In particular, $\pi-\nu$ for $x<0$ and $y\geq 0$, $2\pi-\nu$ for $x\geq 0$ and $y<0$, and $\pi+\nu$ for $x,y<0$ with $\nu$ given by \eqref{nuEll}. 
In figure \ref{fig3}, level sets of cylindrical coordinates $\log{r}$, $\phi$ and elliptic cylindrical coordinates $\mu$, $\nu$ as defined in \eqref{EllC}
are shown together with the respective gradient fields. 
By substituting \eqref{muEll} into the expression \eqref{T}, one obtains elliptic toroidal surfaces.
In the following sections, we will see how to embed symmetric and quasisymmetric magnetic fields within such domains.

\section{Symmetric and self-quasisymmetric magnetic fields in asymmetric tori}
\begin{figure}[h]
\hspace*{-0cm}\centering
\includegraphics[scale=0.46]{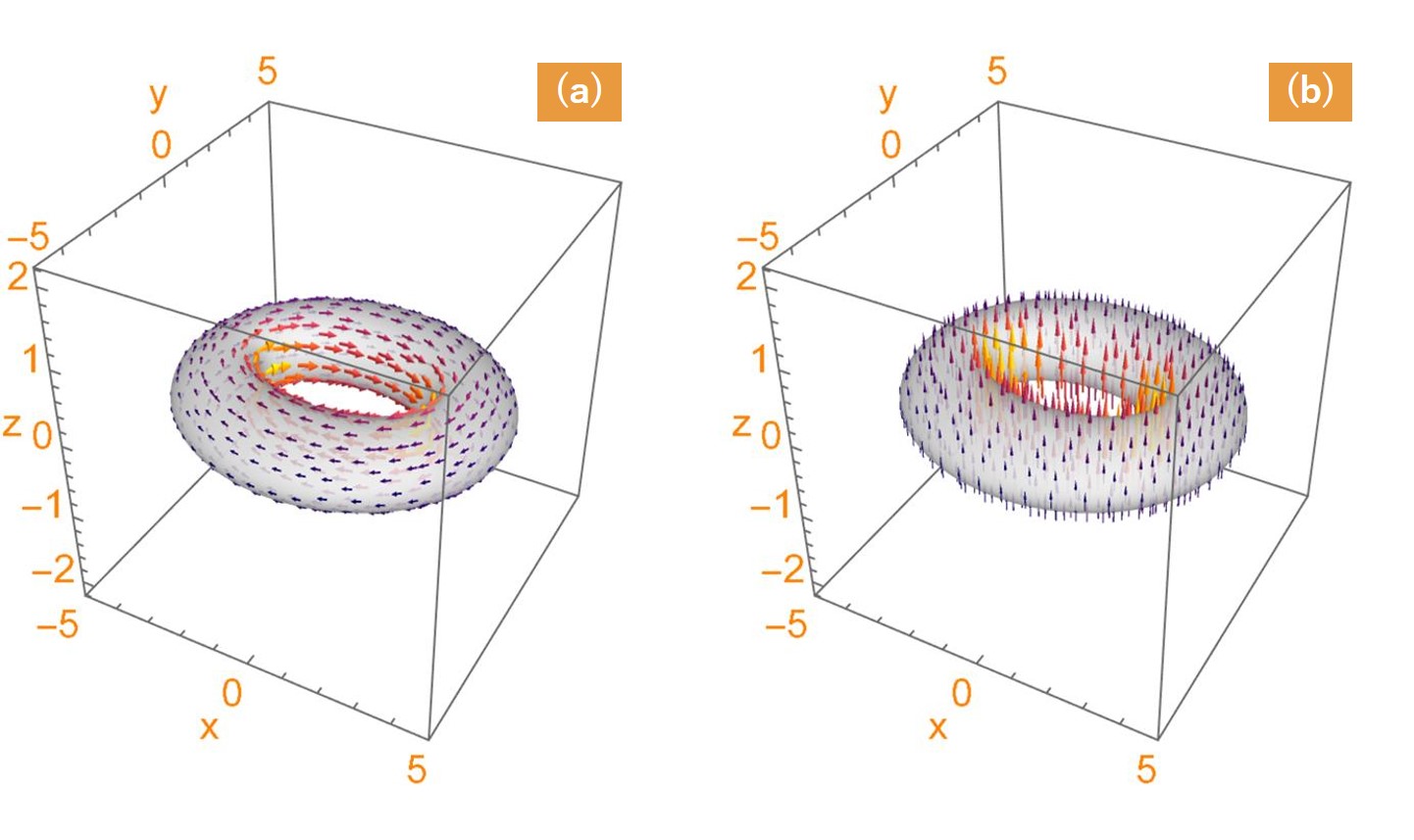}
\caption{\footnotesize (a) Plot of the magnetic field $\bol{B}=-e^{-\mu}\nabla\nu$ with translational symmetry $\bol{u}=\nabla z$ over the surface $\Psi_{\rm el}=0.1$. Here, $\mu$ and $\nu$ are elliptic cylindrical coordinates as given in \eqref{EllC}. This magnetic field is a symmetric solution of anisotropic magnetohydrodynamics \eqref{aMHD2} with pressure fields \eqref{PP} in an asymmetric domain whose boundary is a level set of $\Psi_{\rm el}$. (b) Plot of the corresponding current $\cu\bol{B}=e^{-\mu}\abs{\nabla\mu}^2\nabla z$. In these plots, $\mu_0=1$ and $a=2$.}
\label{fig4}
\end{figure}

Let $\Omega$ denote the volume enclosed by the elliptic torus
\begin{equation}
\Psi_{\rm el}=\frac{1}{2}\left[\lr{\mu-\mu_0}^2+z^2\right],\label{elT}
\end{equation}
with $\mu$ given by \eqref{muEll}. Again, the surface $\Psi_{\rm el}$ does not possess continuous Euclidean symmetries
because the equation $\mf{L}_{\bol{u}}\Psi_{\rm el}=0$ with $\bol{u}=\bol{a}+\bol{b}\cp\bol{x}$ does not have nontrivial solutions $\bol{u}\neq\bol{0}$. 
In particular, note that horizontal cuts of $\Psi_{\rm el}$ are ellipsoidal. 
As discussed in the previous section, in the absence of quasisymmetry  
and boundary conditions on the pressure fields, anisotropic magnetohydrodynamic equilibria within $\Omega$ can be obtained by
enforcing \eqref{PP} and by  setting
$\bol{B}=\nabla\Psi_{\rm el}\cp\nabla\Phi$ with $\Phi$ an arbitrary function whose gradient field is linearly independent of $\nabla\Psi_{\rm el}$. 
Observing that $\nabla\Psi_{\rm el}\cdot\nabla\nu=0$, a family of solutions with a more familiar form is
\begin{equation}
\bol{B}=\nabla\Psi_{\rm el}\cp\nabla\nu+\lambda\nabla\nu,\label{elB}
\end{equation}
with $\lambda=\lambda\lr{\mu,z}$.  
Next, we seek for a translationally symmetric solution such that the direction of symmetry is $\bol{u}=\nabla z$.
For simplicity, assume that 
$\bol{B}=\lambda\nabla\nu$. 
Evidently, $\nabla\cdot\bol{B}=\nabla\cdot\bol{u}=0$. Hence, equation \eqref{aMHD2} is satisfied provided that
\begin{equation}
\lambda\nabla \nu\cp\nabla z=\nabla\zeta,~~~~\frac{\p}{\p z}\lambda^2\abs{\nabla\nu}^2=0.
\end{equation} 
Since by construction $\nabla \nu\cp\nabla z=\nabla\mu$ and $\p\abs{\nabla\nu}/\p z=0$, these conditions are identically satisfied for any $\lambda=\lambda\lr{\mu}$. 
Then, $\zeta=\int\lambda d\mu$.  
This example shows that a symmetric solution can be fitted within an asymmetric domain.
In figure \ref{fig4} a plot of a symmetric magnetic field constructed as described above is given.     

Let us now consider the existence of self-quasisymmetric magnetic fields \eqref{SQS} within $\Omega$ 
such that the perpendicular pressure $P_{\perp}$ is constant on the boundary. 
Such solution can be obtained by solving \eqref{Pb} with $\Psi_{\rm T}=\Psi_{\rm el}$ for
the Clebsch potential $\Phi$. 
Assuming that $B\neq 0$, it is convenient to 
introduce the quantity $\Phi'=\Phi/\sqrt{E\lr{\Psi_{\rm el}}}$ so that
the equation to be solved becomes
\begin{equation}
\abs{\nabla\Psi_{\rm el}\cp\nabla\Phi'}^2=1~~~~{\rm in}~~\Omega.
\end{equation}
Substituting the expression for $\Psi_{\rm el}$, one arrives at 
the following nonlinear first-order partial differential equation
\begin{equation}
\lr{\frac{\p\Phi'}{\p\nu}}^2\left[z^2+\frac{\lr{\mu-\mu_0}^2}{a^2\delta}\right]+\left[z\frac{\p\Phi'}{\p\mu}-\lr{\mu-\mu_0}\frac{\p\Phi'}{\p z}\right]^2=a^2\delta~~~~{\rm in}~~\Omega.\label{Phi1}
\end{equation}
To further simplify the equation, it is convenient to introduce the change of variables $\lr{\mu,z}\rightarrow\lr{\rho,\theta}$ defined by
\begin{equation}
\mu-\mu_0=\rho\cos\theta,~~~~z=\rho\sin\theta.
\end{equation}
Note that $\rho^2=2\Psi_{\rm el}$. 
Then, equation \eqref{Phi1} reads as
\begin{equation}
\lr{\frac{\p\Phi'}{\p\nu}}^2\rho^2\left(\sin^2\theta+\frac{\cos^2\theta}{a^2\delta}\right)+\lr{\frac{\p\Phi'}{\p\theta}}^2=a^2\delta~~~~{\rm in}~~\Omega,\label{Phi2}
\end{equation}
where $\delta=\sin^2\nu+\sinh^2\left(\mu_0+\rho\cos\theta\right)$. 
For given $\rho$, \eqref{Phi2} is a nonlinear first-order partial differential equation
with two independent variables quadratic in the derivatives on a flux surface $\Psi_{\rm el}$. Indeed, it has the form
\begin{equation}
F\lr{p,q,\nu,\theta}=f p^2+q^2-g=0,
\end{equation}
where
\begin{equation}
p=\frac{\p\Phi'}{\p\nu},~~~~q=\frac{\p\Phi'}{\p\theta},~~~~f=\rho^2\left(\sin^2\theta+\frac{\cos^2\theta}{a^2\delta}\right),~~~~g=a^2\delta.
\end{equation}
A local solution of \eqref{Phi2} can be found by considering the Cauchy problem
for the characteristic system of ordinary differential equations associated with \eqref{Phi2} (see e.g. \cite{Polyanin}). 
To this end, denote with $\xi$ a parameter and assign initial conditions as below:
\begin{equation}
\nu=\nu_0\lr{\xi},~~~~\theta=\theta_0\lr{\xi},~~~~\Phi'=\Phi'_0\lr{\xi},~~~~p=p_0\lr{\xi},~~~~q=q_0\lr{\xi},
\end{equation} 
where the functions $\nu_0$, $\theta_0$, $\Phi'_0$, $p_0$, and $q_0$ are determined so that
the following non-degeneracy, compatibility, and transversality conditions are fulfilled:
\begin{subequations}
\begin{align}
&\lr{\frac{\p F}{\p p}}^2+\lr{\frac{\p F}{\p q}}^2\neq 0,\label{c1}\\
&\lr{\frac{d\nu_0}{d\xi}}^2+\lr{\frac{d\theta_0}{d\xi}}^2\neq 0,\label{c2}\\
&F\lr{p_0,q_0,\nu_0,\theta_0}=0,\label{c3}\\
&\frac{d\Phi'_0}{d\xi}=p_0\frac{d\nu_0}{d\xi}+q_0\frac{d\theta_0}{d\xi},\label{c4}\\
&\frac{\p F}{\p p}\frac{d\theta_0}{d\xi}-\frac{\p F}{\p q}\frac{d\nu_0}{d\xi}\neq 0.\label{c5}
\end{align}
\end{subequations}
Equation \eqref{c1} is identically satisfied as long as $f$, $p$, and/or $q$ are different from zero. Indeed, 
\begin{equation}
\lr{\frac{\p F}{\p p}}^2+\lr{\frac{\p F}{\p q}}^2=4f^2p^2+4q^2\neq 0.
\end{equation}
Equation \eqref{c3} can be solved for $q_0$. We have
\begin{equation}
q_0^2=g_0-f_0p^2_0.
\end{equation}
Here, $g_0=g\lr{\nu_0,\theta_0}$ and $f_0=f\lr{\nu_0,\theta_0}$. 
The other equations \eqref{c2}, \eqref{c4}, and \eqref{c5} can be satisfied, for example, by demanding that $\nu_0=\xi$, $\theta_0=c_{\theta}\in\mathbb{R}$,
and $p_0=0$. 
In this case \eqref{c2}, \eqref{c4}, and \eqref{c5} become
\begin{subequations}
\begin{align}
&\lr{\frac{d\nu_0}{d\xi}}^2+\lr{\frac{d\theta_0}{d\xi}}^2=1\neq 0,\label{c2b}\\
&\frac{d\Phi'_0}{d\xi}=p_0\frac{d\nu_0}{d\xi}+q_0\frac{d\theta_0}{d\xi}=0,\label{c4b}\\
&\frac{\p F}{\p p}\frac{d\theta_0}{d\xi}-\frac{\p F}{\p q}\frac{d\nu_0}{d\xi}=-2q_0\neq 0.\label{c5b}
\end{align}
\end{subequations}
Hence, a set of consistent initial conditions is
\begin{equation}
\nu_0=\xi,~~~~\theta_0=c_\theta,~~~~\Phi'_0=c_{\Phi'},~~~~p_0=0,~~~~q_0=\sqrt{g_0}.\label{ic}
\end{equation}
with $c_{\Phi'}\in\mathbb{R}$.
Then, solving the characteristic Lagrange-Charpit system
\begin{equation}
\frac{d\nu}{\frac{\p F}{\p p}}=\frac{d\theta}{\frac{\p F}{\p q}}=\frac{d\Phi'}{p\frac{\p F}{\p p}+q\frac{\p F}{\p q}}=-\frac{dp}{\frac{\p F}{\p\nu}}=-\frac{dq}{\frac{\p F}{\p\theta}}=d\tau,
\end{equation}
one obtains a parametric solution $\nu=\nu\lr{\xi,\tau}$, $\theta\lr{\xi,\tau}$, and $\Phi'\lr{\xi,\tau}$ in a neighborhood of the line defined by the initial conditions \eqref{ic}.
This result implies that, for a given function $E\lr{\Psi_{\rm el}}>0$, a self-quasisymmetric solution $\bol{B}=\sqrt{E\lr{\Psi_{\rm el}}}\nabla\Psi_{\rm el}\cp\nabla\Phi'$ of anisotropic magnetohydrodynamics with constant perpendicular pressure on the boundary  
can be constructed in a neighborhood $U\subset\Omega$ such that the boundary condition $\bol{B}\cdot\bol{n}=0$ is satisfied on $\p U\cap\p\Omega\neq\emptyset$. 
As shown by the example above, in the context of anisotropic magnetohydrodynamics the existence of quasisymmetric solutions 
ultimately translates into the global consistency of locally obtained solutions, i.e. whether the solution defined in $U$ can be
prolonged to cover the whole $\Omega$.  
We therefore conjecture that rigorous results concerning the existence of quasisymmetric vector fields
cannot be obtained unless the issue of global extension of solutions of first-order nonlinear partial differential equations
is carefully addressed.    
Finally, we observe that the results discussed in the present section apply to
asymmetric tori generated trough general orthogonal harmonic coordinates, 
i.e. they are not restricted to elliptic geometry.

\section{Quasisymmetric magnetic fields in asymmetric tori}
The purpose of the present section is to derive quasisymmetric magnetic fields. 
We will see that these vector fields can be interpreted as local solutions of anisotropic magnetohydrodynamics \eqref{aMHD2} 
in an asymmetric torus with perpendicular and parallel pressure fields given by \eqref{PP}. 
In particular, these solutions hold in a neighborhood $U\subset\Omega$ with $\p U\cap\p\Omega\neq\emptyset$, they do not have continuous Euclidean symmetries, no boundary conditions are imposed on the pressure fields, and $\bol{B}\cp\bol{u}=\nabla\zeta\neq\bol{0}$ (the fields are not self-quasisymmetric).
To achieve this goal, we solve equation \eqref{aMHD3} in a neighborhood $U$ of a toroidal domain $\Omega$ 
whose boundary $\p\Omega$ corresponds to a level set of a function $\Psi_{\rm T}$. 

First, we prescribe the toroidal domain $\Omega$. 
Instead of breaking axial symmetry by modifying the distance function $r$ with a new variable $\mu$, 
it is convenient to deform the torus by introducing a non-zero $h$ in ${\Psi_{\rm T}}$ of \eqref{T}. 
To simplify calculations, we postulate that $h=h\lr{r,\phi}$.
We thus have
\begin{equation}
\Psi_{\rm T}=\frac{1}{2}\left[\lr{r-r_0}^2+\lr{z-h}^2\right].\label{PsiTQS1}
\end{equation} 
For the magnetic field to possess nested flux surfaces, we further demand that
\begin{equation}
\bol{B}=f\lr{r,z-h}\nabla r\cp\nabla\lr{z-h}.\label{BQS1}
\end{equation}
Here, $f$ is a function of $r$ and $z-h$ to be determined together with $h$ by imposing the quasisymmetry condition.
Notice that $\bol{B}\cdot\nabla\Psi_{\rm T}=0$ implying $\bol{B}\cdot\bol{n}=0$ on $\p\Omega$. The magnetic field also satisfies $\nabla\cdot\bol{B}=0$.
Hence, the remaining equations in \eqref{aMHD2} are
\begin{equation}
\bol{B}\cp\bol{u}=\frac{d g}{d\zeta}\nabla\zeta,~~~~\bol{u}\cdot\nabla B=0,~~~~\nabla\cdot\bol{u}=0~~~~{\rm in}~~\Omega.\label{Bu}
\end{equation} 
In the equation above we have replaced $\zeta$ with a function of $\zeta$, $g=g\lr{\zeta}$, for later convenience. 
Assuming the Clebsch representation $\bol{u}=\nabla\zeta\cp\nabla\rho$ with $\zeta=z-h$, 
from \eqref{Bu} we thus arrive at \eqref{aMHD3}, which now reads
as a system of equations that must be solved for the Clebsch potential $\rho$ and the displacement $h$,
\begin{subequations}
\begin{align}
&\nabla r\cp\nabla\lr{z-h}\cdot\nabla\rho=1,\label{Bu1}\\
&\nabla\lr{z-h}\cp\nabla\rho\cdot\nabla\left[\frac{1}{r^2}\lr{\frac{\p h}{\p\phi}}^2\right]=0.\label{Bu2}
\end{align}
\end{subequations}
Observe that in deriving \eqref{Bu1} and \eqref{Bu2} we set $dg/d\zeta=f$ and
eliminated the dependence of $f$ on $r$.  
Next, 
define the quantity
\begin{equation}
\eta=\frac{1}{r}\frac{\p h}{\p\phi}.
\end{equation}
Then, assuming that $\nabla\lr{z-h}$ and $\nabla\lr{r^{-1}\p h/\p\phi}^2$ are linearly independent, equation \eqref{Bu2} implies that
\begin{equation}
\rho=\rho\lr{z-h,\eta}.\label{Bu22}
\end{equation}
Substituting \eqref{Bu22} into \eqref{Bu1} we arrive at
\begin{equation}
\frac{\p\rho}{\p\eta}\nabla r\cp\nabla\lr{z-h}\cdot\nabla\eta=1,\label{eta1}
\end{equation}
which, after some manipulations, becomes
\begin{equation}
-\frac{\p\rho}{\p\eta}\frac{1}{r^2}\frac{\p^2 h}{\p\phi^2}=1,\label{h}
\end{equation}
Since $h$ does not depend on $z$, we put $\rho=\rho\lr{\eta}$. 
Then, for given $\rho\lr{\eta}$, equation \eqref{h} must be solved for $h$. 
Now recall that $\p h/\p\phi=r\eta$. Therefore, equation \eqref{eta1} is equivalent to
\begin{equation}
\frac{\p\rho}{\p\phi}=-r.
\end{equation}
Assuming that the function $\rho=\rho\lr{\eta}$ can be inverted with inverse $\rho^{-1}$ and integrating with respect to $\phi$ gives
\begin{equation}
\frac{\p h}{\p\phi}=r\rho^{-1}\lr{\alpha-r\phi},\label{h2}
\end{equation} 
where $\alpha=\alpha\lr{r}$ is an arbitrary integration factor.
Once $\rho^{-1}$ is assigned, a further integration with respect to $\phi$ gives the desired solution. 
However, due to the radial dependence of $\alpha-r\phi$, the quantity $\p h/\p\phi$ 
is not periodic in $\phi$. Therefore, the components of the corresponding quasisymmetric magnetic field \eqref{BQS1} and the
flux function \eqref{PsiTQS1}  
become multivalued functions, and the obtained quasisymmetric solution is well-defined only locally. 
In particular, the solution is defined in some region $U\subset\Omega$, and the irregularities 
of the bounding surface $\p\Omega$ arising from the multivalued nature of $\Psi_{\rm T}$ 
can be removed by repairing (smoothing) the flux function $\Psi_{\rm T}$ outside $U$.     
To see this, set $\rho\lr{\eta}=k^{-1}\eta$, with $k>0$ a real constant, and choose $f=\sin\lr{z-h}$.
Then, we have
\begin{equation}
\frac{\p h}{\p\phi}=k r\lr{\alpha-r\phi},~~~~h=k r\phi\lr{\alpha-\frac{1}{2}r\phi}+\beta,
\end{equation}
with $\beta=\beta\lr{r}$ an arbitrary integration factor, and also
\begin{subequations}
\begin{align}
&\bol{B}=-\sin\left[z-\beta-k r\phi\lr{\alpha-\frac{1}{2}r\phi}\right]\left[r\nabla\phi+k\lr{\alpha-r\phi}\nabla z\right],\\
&\bol{u}=r\lr{\frac{\p\alpha}{\p r}-\phi}\nabla\phi+\nabla r+\frac{\p}{\p r}\left(k\frac{\alpha^2}{2}+\beta\right)\nabla z.
\end{align}
\end{subequations}  
For the component of the magnetic field to be single-valued, one must determine $\alpha$ and $\beta$ so that
\begin{equation}
\bol{B}\rvert_{\phi=0}=\bol{B}\rvert_{\phi=2\pi},
\end{equation}
which implies
\begin{subequations}
\begin{align}
&z-\beta=z-\beta-2\pi k r\lr{\alpha-\pi r}+2\pi m,~~~~m\in\mathbb{Z}\\
&k\alpha=k\lr{\alpha-2\pi r}.
\end{align}
\end{subequations}
However, excluding the case $k=0$ which corresponds to axial symmetry, 
these conditions cannot be satisfied. Hence, the local nature of the solution.

\begin{figure}[h]
\hspace*{-0cm}\centering
\includegraphics[scale=0.5]{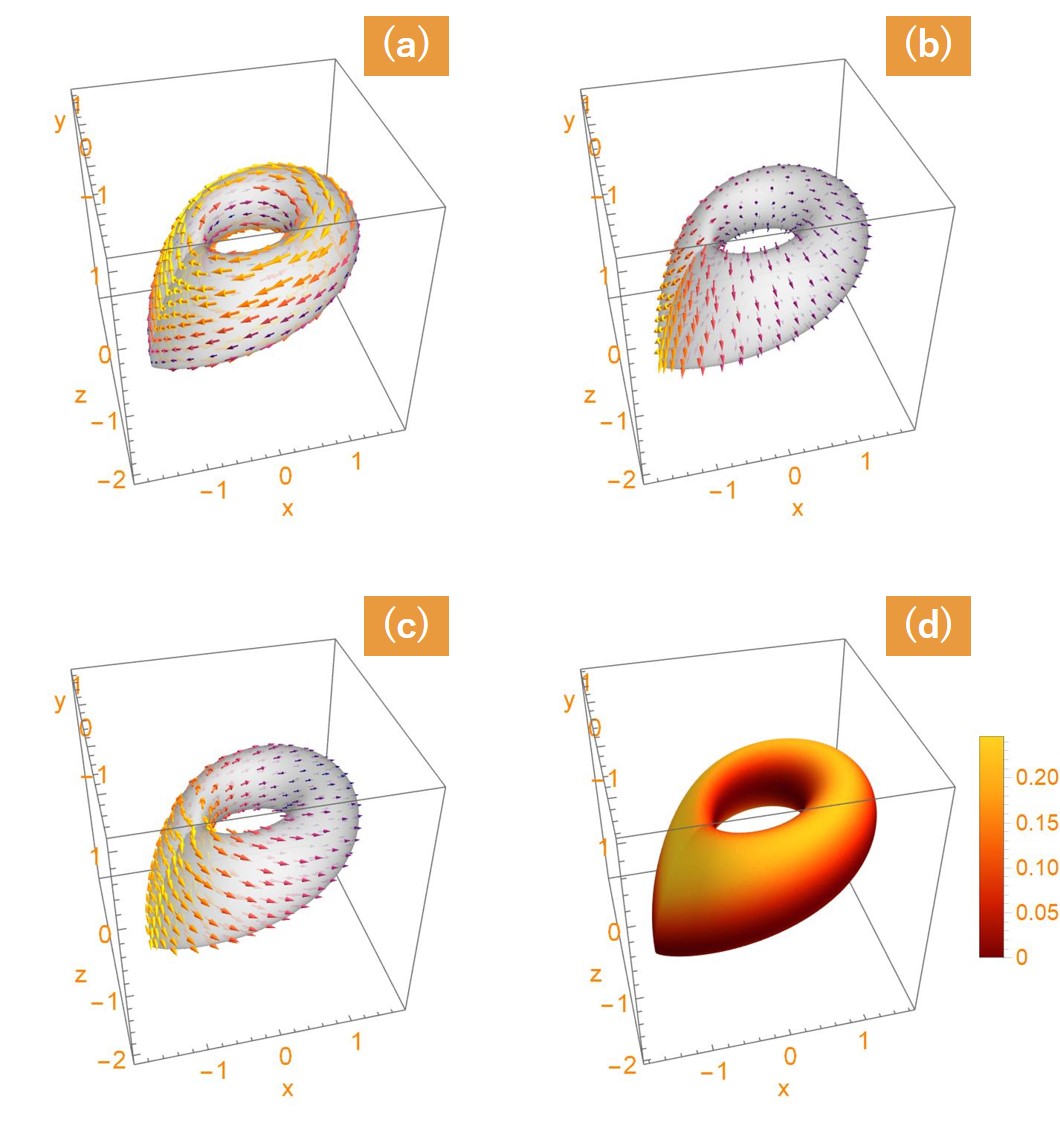}
\caption{\footnotesize (a) Plot of the quasisymmetric magnetic field $\bol{B}$ of \eqref{BQS2} on a level set of \eqref{PsiTQS1}. (b) Plot of the current $\nabla\cp\bol{B}$. (c) Plot of the quasisymmetry $\bol{u}$. (d) Plot of the field stregth $B^2$. Observe that $\bol{B}$, $\cu\bol{B}$, $\bol{u}$, and $\Psi_{\rm T}$ exhibit a discontinuity at $\phi=0$ due to the multivalued nature of $\phi$. In these plots, $k=0.18$.}
\label{fig5}
\end{figure}

As an example of quasisymmetric configuration, consider the case $\alpha=\beta=0$. Then, one can verify the quasisymmetry of the derived solution as below:
\begin{subequations}\label{BQSEx}
\begin{align}
&\bol{B}=-\sin\left(z+\frac{1}{2}k r^2\phi^2\right)\left(r\nabla\phi-k r\phi\nabla z\right),\label{BQS2}\\
&\nabla\cdot\bol{B}=0,\\
&B^2=\sin^2\lr{z+\frac{1}{2}k r^2\phi^2}\lr{1+k^2r^2\phi^2},\\
&\bol{u}=\nabla r-r\phi\nabla\phi,\\
&\nabla\cdot\bol{u}=0,\\
&\bol{B}\cp\bol{u}=-\nabla\cos\lr{z+\frac{1}{2}k r^2\phi^2},\\
&\bol{u}\cdot\nabla B^2=0,\\
&\bol{B}\cdot\nabla\Psi_{\rm T}=0.
\end{align}
\end{subequations}
We also remark that this solution does not have continuous Euclidean symmetries
because the equation $\mf{L}_{\bol{u}}\bol{B}=\bol{0}$ with $\bol{u}=\bol{a}+\bol{b}\cp\bol{x}$ does not have solutions with $\bol{a},\bol{b}\neq\bol{0}$. 
A plot of the obtained magnetic field and related quantities is given in figure \ref{fig5}.  

The quasisymmetric magnetic field of equation \eqref{BQSEx} is such that the quasisymmetry $\bol{u}$ is not tangential to flux surfaces 
$\Psi_{\rm T}$ since $\bol{u}\cdot\nabla\Psi_{\rm T}\neq 0$. However, in the design of a stellarator it is customary to demand that $\bol{u}$ lies on flux surfaces
to ensure that the conserved momentum associated with guiding center motion is a combination of the poloidal and toroidal
momenta. This requirement amounts to identifying $\zeta$ with the flux function $\Psi_{\rm T}$ in the quasisymmetry equations \eqref{QS2}. 
Such configuration can be achieved by repeating the procedure discussed above while enforcing the Clebsch representations
\begin{equation}
\bol{B}=f\lr{\Psi_{\rm T}}\nabla r\cp\nabla\lr{z-h},~~~~\bol{u}=\nabla\Psi_{\rm T}\cp\nabla\rho,
\end{equation}  
where $f$ is an arbitrary function of $\Psi_{\rm T}$. In this case, system \eqref{aMHD3} reduces to
\begin{subequations}
\begin{align}
&\nabla r\cp\nabla\lr{z-h}\cdot\nabla\rho=1,\\
&\nabla\Psi_{\rm T}\cp\nabla\rho\cdot\nabla\left[\frac{1}{r^2}\lr{\frac{\p h}{\p\phi}}^2\right]=0.
\end{align}
\end{subequations}
Here, we set $dg/d\zeta=dg/d\Psi_{\rm T}=f$. 
Again, defining $\eta=r^{-1}\p h/\p\phi$ and requiring that $\rho=\rho\lr{\eta}$, one arrives at equation \eqref{h} with solution \eqref{h2}.
For example, if $\rho=k^{-1}\eta$ and $\alpha=\beta=0$ one obtains the family of quasisymmetry magnetic fields
\begin{subequations}\label{BQSEx2}
\begin{align}
&\bol{B}=-f\left(r\nabla\phi-k r\phi\nabla z\right),\\
&\nabla\cdot\bol{B}=0,\\
&B^2=f^2\lr{1+k^2r^2\phi^2},\\
&\bol{u}=-\lr{r-r_0}\nabla z+\lr{z+\frac{1}{2}k r^2\phi^2}\lr{\nabla r-r\phi\nabla\phi},\\
&\nabla\cdot\bol{u}=0,\\
&\bol{B}\cp\bol{u}=f\nabla\Psi_{\rm T}=\nabla g,\\
&\bol{u}\cdot\nabla B^2=0,\\
&\bol{B}\cdot\nabla\Psi_{\rm T}=0,\\
&\bol{u}\cdot\nabla\Psi_{\rm T}=0.
\end{align}
\end{subequations}
Observe that both $\bol{B}$ and $\bol{u}$ are tangential to flux surfaces. 

The construction of quasisymmetric magnetic fields presented above can be generalized to the case of  
surfaces defined through harmonic orthogonal coordinates $\lr{\mu,\nu,z}$. 
Setting $h=h\lr{\mu,\nu}$, the relevant Clebsch parametrization is  
\begin{equation}
\Psi_{\rm T}=\frac{1}{2}\left[\lr{\mu-\mu_0}^2+\lr{z-h}^2\right],~~~~\bol{B}=f\lr{\Psi_{\rm T}}\nabla\mu\cp\nabla \lr{z-h},~~~~\bol{u}=\nabla\Psi_{\rm T}\cp\nabla\rho.\label{QSF}
\end{equation}
Then, putting $dg/d\Psi_{\rm T}=f$, the quasisymmetry equations \eqref{aMHD3} now read
\begin{subequations}
\begin{align}
&\nabla\Psi_{\rm T}\cp\nabla\rho\cdot\nabla\left\{\abs{\nabla\mu}^2\left[1+\abs{\nabla\mu}^2\lr{\frac{\p h}{\p\nu}}^2\right] \right\}=0,\\
&\nabla\mu\cp\nabla\lr{z-h}\cdot\nabla\rho=1.
\end{align}
\end{subequations}
Again, requiring that $\rho=\rho\lr{\eta}$ with
\begin{equation}
\eta=\abs{\nabla\mu}^2\left[1+\abs{\nabla\mu}^2\lr{\frac{\p h}{\p\nu}}^2\right],
\end{equation} 
we arrive at one equation which determines the displacement $h$:
\begin{equation}
-\frac{d\rho}{d\eta}\frac{\p\eta}{\p\nu}\abs{\nabla\mu}^2=1.
\end{equation}
For a given $\rho$, the corresponding solution $h$ of the equation above assigns a family of quasisymmetric magnetic fields parametrized by the function $f$ 
as defined in equation \eqref{QSF}. 
Nonetheless, as in the previous cases, the globality of the solution is not guaranteed.

\section{Concluding remarks}
In this work, we studied the existence of quasisymmetric magnetic fields in asymmetric toroidal domains. 
If the perpendicular and parallel pressure fields are chosen as in equation \eqref{PP}, 
this problem is equivalent to finding quasisymmetric magnetohydrodynamic equilibria within the framework of anisotropic magentohydrodynamics. 
In particular, constructing a quasisymmetric configuration amounts to solving system \eqref{aMHD2} for the vector fields $\bol{B}$ and $\bol{u}$. 
This system can be written in terms of two coupled nonlinear first-order partial differential equations \eqref{aMHD3}, 
or a single nonlinear second-order partial differential equation \eqref{GSQS}.  
By using harmonic orthogonal coordinates, we first devised a method to obtain symmetric solutions of \eqref{aMHD3} in asymmetric toroidal domains.  
Then, we constructed regular local self-quasisymmetric and quasisymmetric magnetic fields in asymmetric tori 
by solving \eqref{aMHD3} through a Clebsch parametrization of the solution tailored on the flux function associated with the toroidal boundary.  
These solutions are local in the sense that they solve \eqref{aMHD3} in a region $U\subset\Omega$
and satisfy boundary conditions on a portion of the boundary, $\p U\cap\p\Omega\neq\emptyset$.  

The obtained results highlight two aspects that characterize quasisymmetric vector fields.  
On one hand, the symmetry of the boundary should be treated separately from the symmetry of the solution, since
we have shown that a symmetric magnetic field can be fitted within an asymmetric domain. 
On the other hand, due to the first-order nature of the governing equations, 
the mathematical challenge posed by quasisymmetry can be ascribed to the 
local nature of the solutions obtained by integrating the characteristic system of ordinary
differential equations. 
Finally, we note that the possibility that the obstruction 
encountered in the derivation of global solutions is intrinsic to three-dimensional equilibria cannot be ruled out, 
in the sense that the absence of Euclidean isometry may prevent the existence of such fields.
In particular, it appears that the existence of global quasisymmetric fields
is contingent upon the availability of curvilinear coordinates whose 
metric coefficients exhibit invariance properties analogous to those
satisfied by cylindrical coordinates, and specifically by the toroidal angle $\phi$ which obeys $\p\abs{\nabla\phi}/\p\phi=0$.  




\section*{Acknowledgment}
The research of NS was supported by JSPS KAKENHI Grant No. 21K13851 and 17H01177.

\section*{Data Availability}

The data that support the findings of this study are available from
the corresponding author upon reasonable request.

\end{document}